\documentclass{amsart}

\usepackage{fancyhdr}
\usepackage{lastpage}
\usepackage{stmaryrd,yhmath}
\usepackage{tikz}
\usepackage[all]{xy}

\newcommand{\bdis}{\begin{displaymath}}
\newcommand{\edis}{\end{displaymath}}
\newcommand{\be}{\begin{equation}}
\newcommand{\ee}{\end{equation}}
\newcommand{\mbb}{\mathbb}
\newcommand{\mcal}{\mathcal}

\theoremstyle{definition}

\theoremstyle{remark}
\newtheorem{remark}[]{Remark}

\newtheorem*{mydef1}{{\bf Theorem}}

\newtheorem*{mydef41}{{\bf Corollary 1}}

\newtheorem*{mydef42}{{\bf Corollary 2}}

\newtheorem*{mydef43}{{\bf Corollary 3}}

\newtheorem*{mydef44}{{\bf Corollary 4}}

\newtheorem*{mydef51}{{\bf Lemma 1}}

\newtheorem*{mydef52}{{\bf Lemma 2}}

\newtheorem*{mydef53}{{\bf Lemma 3}}

\numberwithin{equation}{section}

\begin{document}

\title{Jacob's ladders, nonlinear interactions between $\zeta$-oscillating systems and corresponding constraints}

\author{Jan Moser}

\address{Department of Mathematical Analysis and Numerical Mathematics, Comenius University, Mlynska Dolina M105, 842 48 Bratislava, SLOVAKIA}

\email{jan.mozer@fmph.uniba.sk}


\begin{abstract}
In this paper we introduce new class of nonlinear interactions of $\zeta$-oscillating systems. The main formula is generated by corresponding subset
of the set of trigonometric functions. Next, the main formula generates certain set of two-parts forms. For this set the following holds true: the cube
of two-part form is asymptotically equal to other two-part form -- short functional algebra.
\end{abstract}
\maketitle

\section{Introduction}

\subsection{}

Let us remind that we have obtained in our paper \cite{8} certain class of formulae linear in different variables of the following type
\be \label{1.1}
\prod_{r=1}^k\left|\frac{\zeta\left(\frac 12+i\alpha_r\right)}{\zeta\left(\frac 12+i\beta_r\right)}\right|^2,\ k=1,\dots,k_0,
\ee
i.e. linear in different oscillating systems. Namely, see formulae in \cite{8}:
\be \label{1.2}
(7.1),\ (7.3),\ (7.4),\ (8.2), \ (8.4), \ (8.5).
\ee
It is, for example, the formula (7.1):
\be \label{1.3}
\begin{split}
& \cos^2(\alpha_0^{2,2})\prod_{r=1}^{k_2}\left|\frac{\zeta\left(\frac 12+i\alpha_r^{2,2}\right)}{\zeta\left(\frac 12+i\beta_r^2\right)}\right|^2+ \\
& + \sin^2(\alpha_0^{1,1})\prod_{r=1}^{k_1}\left|\frac{\zeta\left(\frac 12+i\alpha_r^{1,1}\right)}{\zeta\left(\frac 12+i\beta_r^1\right)}\right|^2\sim 1, \\
& L\to\infty,\ 1\leq k_1,k_2\leq k_0,\ L,k_0\in\mbb{N},
\end{split}
\ee
where
\bdis
\begin{split}
& \alpha_0^{1,1}=\alpha_0^{1}(U,\mu,L,k_1;\sin^2t), \dots \\
& \beta_r^1=\beta_r(U,\mu,L,k_1), \\
& \alpha_0^{2,2}=\alpha_0^{2}(U,\mu,L,k_2;\cos^2t), \dots \\
& \beta_r^2=\beta_r(U,\mu,L,k_2) .
\end{split}
\edis

\begin{remark}
We have called the formula (\ref{1.2}) as the $\zeta$-analogue of the elementary trigonometric identity
\bdis
\cos^2t+\sin^2t=1.
\edis
\end{remark}

Next, we have introduced in our paper \cite{8} the following notions in connection with the formula (\ref{1.2}):
\begin{itemize}
\item[(a)] functionally depending $\zeta$-oscillating systems and linearly connected $\zeta$-oscillating systems, (see \cite{8}, beginning of section 2.4),
\item[(b)] interaction between corresponding $\zeta$-oscillating systems (see \cite{8}, Definition 4).
\end{itemize}

\begin{remark}
By (a) and (b) we may assume in the context of the paper \cite{8} the following:
\begin{center}
interaction \ = \ linear interaction,
\end{center}
(see, for example, (\ref{1.3}), comp. \cite{9}, Remark 1).
\end{remark}

\begin{remark}
Moreover, we notice that the $\zeta$-oscillating system itself (comp. (\ref{1.1})) is a complicated nonlinear system (comp. \cite{8}, (1.7), i.e. the spectral
form of the Riemann-Siegel formula).
\end{remark}

\subsection{}

Next, let us remind the following notions we have introduced (see \cite{1} -- \cite{7}) within the theory of the Riemann zeta-function:
\begin{itemize}
\item[(A)] Jacob's ladders, (see \cite{1}, comp. \cite{2}),
\item[(B)] $\zeta$-oscillating system, (see \cite{7}, (1.1)),
\item[(C)] factorization formula, (see \cite{5}, (4.3) -- (4.18), comp. \cite{7}, (2.1) -- (2.7)),
\item[(D)] metamorphosis of the $\zeta$-oscillating systems:
\begin{itemize}
\item[(a)] first, the notion of metamorphoses of an oscillating multiform \cite{4},
\item[(b)] after that, the notion of metamorphoses of a quotient of two oscillating multiform, \cite{5},
\end{itemize}
\item[(E)] $\mcal{Z}_{\zeta,Q^2}$-transformation (or device), \cite{7}.
\end{itemize}

\subsection{}

In this paper, we shall present certain class of nonlinear interaction formulae, namely:
\begin{itemize}
\item[(a)] containing nonlinearities in variables of kind (\ref{1.1}), (comp. Remark 2),
\item[(b)] describing interaction between corresponding $\zeta$-oscillating systems, i.e. every of these is functionally depending upon
others of these $\zeta$-oscillating systems.
\end{itemize}

The main result is expressed by the following nonlinear formula:
\be \label{1.4}
\begin{split}
& \prod_{r=1}^{k_1}\left|\frac{\zeta\left(\frac 12+i\alpha_r^{1,k_1}\right)}{\zeta\left(\frac 12+i\beta_r^{k_1}\right)}\right|^2\sim \\
& \sim \frac{\cos(\alpha_0^{2,k_2})}{\cos^3(\alpha_0^{1,k_1})}
\prod_{r=1}^{k_2}\left|\frac{\zeta\left(\frac 12+i\alpha_r^{2,k_2}\right)}{\zeta\left(\frac 12+i\beta_r^{k_2}\right)}\right|^2- \\
& - \frac{U^2}{3}\frac{1}{\cos^3(\alpha_0^{1,k_1})\cos^3U}\times \\
& \times
\left\{
\cos^2(\alpha_0^{3,k_3})\prod_{r=1}^{k_3}
\left|\frac{\zeta\left(\frac 12+i\alpha_r^{3,k_3}\right)}{\zeta\left(\frac 12+i\beta_r^{k_3}\right)}\right|^2-
\sin^2(\alpha_0^{4,k_4})\prod_{r=1}^{k_4}
\left|\frac{\zeta\left(\frac 12+i\alpha_r^{4,k_4}\right)}{\zeta\left(\frac 12+i\beta_r^{k_4}\right)}\right|^2
\right\}^3 , \\
& L\to\infty,\ 1\leq k_1,k_2,k_3,k_4\leq k_0,
\end{split}
\ee
(comp. (\ref{1.3}), (\ref{1.4}) and (a) in the beginning of this section).

\section{Lemmas}

\subsection{}

Since
\bdis
\int\cos^3t{\rm d}t=\sin t-\frac 13\sin^3t ,
\edis
then
\bdis
\int_{2\pi L}^{2\pi L+U}\cos^3t{\rm d}t=\sin U-\frac 13\sin^3U,\ L\in\mbb{N},
\edis
and
\bdis
\frac 1U\int_{2\pi L}^{2\pi L+U}\cos^3t{\rm d}t=\frac{\sin U}{U}-\frac{U^2}{3}\left(\frac{\sin U}{U}\right)^3.
\edis
Of course,
\be \label{2.1}
f_1(t)=\cos^3t\in \tilde{C}_0[2\pi L,2\pi L+U],
\ee
(the class of functions $\tilde{C}_0$ has been defined in our paper \cite{9}), where
\bdis
t\in [2\pi L,2\pi L+U],\ U\in\left(\left. 0,\frac{\pi}{2}-\epsilon\right.\right],
\edis
and $\epsilon>0$ is sufficiently small. Consequently, we obtain for generating of the factorization formula by making use our algorithm (see \cite{8}, (3.1) -- (3.11)) the following

\begin{mydef51}
For the function (\ref{2.1}) there are vector-valued functions
\bdis
(\alpha_0^{1,k_1},\alpha_1^{1,k_1},\dots,\alpha_{k_1}^{1,k_1},\beta_1^{k_1},\dots,\beta_{k_1}^{k_1}),\ k_1=1,\dots,k_0,\ k_0\in\mbb{N}
\edis
($k_0$ is arbitrary and fixed) such that the following factorization formula holds true:
\be \label{2.2}
\begin{split}
& \prod_{r=1}^{k_1}\left|\frac{\zeta\left(\frac 12+i\alpha_r^{1,k_1}\right)}{\zeta\left(\frac 12+i\beta_r^{k_1}\right)}\right|^2\sim \\
& \sim \left[\frac{\sin U}{U}-\frac{U^2}{3}\left(\frac{\sin U}{U}\right)^3\right]\frac{1}{\cos^3(\alpha_0^{1,k_1})},\ L\to\infty,
\end{split}
\ee
where
\be \label{2.3}
\begin{split}
& \alpha_r^{1,k_1}=\alpha_r(U,L,k_1;f_1),\ r=0,1,\dots,k_1, \\
& \beta_r^{k_1}=\beta_r(U,L,k_1),\ r=1,\dots,k_1, \\
& 2\pi L<\alpha_0^{1,k_1}<2\pi L+U \ \Rightarrow \ 0<\alpha_0^{1,k_1}-2\pi L<U, \\
& 1\leq k_1\leq k_0.
\end{split}
\ee
\end{mydef51}

\begin{remark}
In the asymptotic formula (\ref{2.2}) the symbol $\sim$ stands for (see \cite{8}, (3.8))
\bdis
=\left\{ 1+\mcal{O}\left(\frac{\ln\ln L}{\ln L}\right)\right\}.
\edis
\end{remark}

\subsection{}

Let
\be \label{2.4}
f_2(t)=\cos t,\ t\in [2\pi L,2\pi L+U],\ U\in\left(\left. 0,\frac{\pi}{2}-\epsilon\right.\right].
\ee
Since
\bdis
f_2(t)\in\tilde{C}_0[2\pi L,2\pi L+U]
\edis
then we obtain by similar way as in the case of Lemma 1 the following.

\begin{mydef52}
For the function (\ref{2.4}) there are vector-valued functions
\bdis
(\alpha_0^{2,k_2},\alpha_1^{2,k_2},\dots,\alpha_{k_2}^{2,k_2},\beta_1^{k_2},\dots,\beta_{k_2}^{k_2}),\ k_2=1,\dots,k_0,\ k_0\in\mbb{N}
\edis
($k_0$ is arbitrary and fixed) such that the following factorization formula holds true:
\be \label{2.5}
\begin{split}
& \prod_{r=1}^{k_2}\left|\frac{\zeta\left(\frac 12+i\alpha_r^{2,k_2}\right)}{\zeta\left(\frac 12+i\beta_r^{k_2}\right)}\right|^2\sim
\frac{\sin U}{U}\frac{1}{\cos(\alpha_0^{2,k_2})},\ L\to\infty,
\end{split}
\ee
where
\be \label{2.6}
\begin{split}
& \alpha_r^{2,k_2}=\alpha_r(U,L,k_2;f_2),\ r=0,1,\dots,k_2, \\
& \beta_r^{k_2}=\beta_r(U,L,k_2),\ r=1,\dots,k_2, \\
& 0<\alpha_0^{2,k_2}-2\pi L<U, \\
& 1\leq k_2\leq k_0.
\end{split}
\ee
\end{mydef52}

\subsection{}

Next, let us remind (see \cite{8}, (8.1) and also (4.2), (4.3), (4.6), (4.7)) the following formula
\be \label{2.7}
\begin{split}
& \frac{\cos^2(\alpha_0^{2,2})}{\cos(2\mu+U)}
\prod_{r=1}^{k_2}\left|\frac{\zeta\left(\frac 12+i\alpha_r^{2,2}\right)}{\zeta\left(\frac 12+i\beta_r^2\right)}\right|^2-
\frac{\sin^2(\alpha_0^{1,1})}{\cos(2\mu+U)}
\prod_{r=1}^{k_1}\left|\frac{\zeta\left(\frac 12+i\alpha_r^{1,1}\right)}{\zeta\left(\frac 12+i\beta_r^1\right)}\right|^2\sim \\
& \sim \frac{\sin U}{U},\ L\to \infty, \\
& t\in [\pi L+\mu, \pi L+\mu +U].
\end{split}
\ee
Now, if we put (in our present context)
\bdis
\begin{split}
& \mu=0, L\rightarrow 2L, \\
& \alpha_r^{2,2} \rightarrow \alpha_r^{3,k_3},\ r=0,1,\dots,k_3, \\
& \beta_r^2 \rightarrow \beta_r^{k_3},\ r=1,\dots,k_3, \\
& \alpha_r^{1,1} \rightarrow \alpha_r^{4,k_4},\ r=0,1,\dots,k_4, \\
& \beta_r^1 \rightarrow \beta_r^{k_4},\ r=1,\dots,k_4,
\end{split}
\edis
then we obtain from (\ref{2.7}) in the case
\bdis
t\in [2\pi L, 2\pi L+U],\ U\in\left(\left. 0,\frac{\pi}{2}-\epsilon\right.\right]
\edis
the following

\begin{mydef53}
\be \label{2.8}
\begin{split}
& \frac{\cos^2(\alpha_0^{3,k_3})}{\cos U}\prod_{r=1}^{k_3}
\left|\frac{\zeta\left(\frac 12+i\alpha_r^{3,k_3}\right)}{\zeta\left(\frac 12+i\beta_r^{k_3}\right)}\right|^2- \\
& - \frac{\sin^2(\alpha_0^{4,k_4})}{\cos U}\prod_{r=1}^{k_4}
\left|\frac{\zeta\left(\frac 12+i\alpha_r^{4,k_4}\right)}{\zeta\left(\frac 12+i\beta_r^{k_4}\right)}\right|^2\sim \frac{\sin U}{U},\ L\to\infty,
\end{split}
\ee
where
\be \label{2.9}
\begin{split}
& f_3(t)=\cos^2t \rightarrow \prod_{r=1}^{k_3}
\left|\frac{\zeta\left(\frac 12+i\alpha_r^{3,k_3}\right)}{\zeta\left(\frac 12+i\beta_r^{k_3}\right)}\right|^2 , \\
& f_4(t)=\sin^2t \rightarrow \prod_{r=1}^{k_4}
\left|\frac{\zeta\left(\frac 12+i\alpha_r^{4,k_4}\right)}{\zeta\left(\frac 12+i\beta_r^{k_4}\right)}\right|^2, \\
& t\in [2\pi L, 2\pi L+U],
\end{split}
\ee
and
\bdis
\begin{split}
& \alpha_0^{3,k_3}=\alpha_0^{2,2}(U,0,2L,k_3;f_3),\dots \\
& \beta_r^{k_3}=\beta_r^{2}(U,0,2L,k_3), \\
& \alpha_0^{4,k_4}=\alpha_0^{1,1}(U,0,2L,k_4;f_4),\dots \\
& \beta_r^1 \rightarrow \beta_r^{k_4},\ r=1,\dots,k_4, \\
& \beta_r^{k_4}=\beta_r^{1}(U,0,2L,k_4).
\end{split}
\edis
\end{mydef53}

\section{Theorem}

\subsection{}

Now, we obtain by making use of (\ref{2.2}), (\ref{2.5}) and (\ref{2.7}) the following nonlinear interaction formula.

\begin{mydef1}
\be \label{3.1}
\begin{split}
& \prod_{r=1}^{k_1}\left|\frac{\zeta\left(\frac 12+i\alpha_r^{1,k_1}\right)}{\zeta\left(\frac 12+i\beta_r^{k_1}\right)}\right|^2\sim \\
& \sim \frac{\cos(\alpha_0^{2,k_2})}{\cos^3(\alpha_0^{1,k_1})}\prod_{r=1}^{k_2}
\left|\frac{\zeta\left(\frac 12+i\alpha_r^{2,k_2}\right)}{\zeta\left(\frac 12+i\beta_r^{k_2}\right)}\right|^2-\frac{U^2}{3}\frac{1}{\cos^2(\alpha_0^{1,k_1})\cos^3U}\times \\
& \times
\left\{
\cos^2(\alpha_0^{3,k_3})\prod_{r=1}^{k_3}\left|\frac{\zeta\left(\frac 12+i\alpha_r^{3,k_3}\right)}{\zeta\left(\frac 12+i\beta_r^{k_3}\right)}\right|^2-
\sin^2(\alpha_0^{4,k_4})\prod_{r=1}^{k_4}
\left|\frac{\zeta\left(\frac 12+i\alpha_r^{4,k_4}\right)}{\zeta\left(\frac 12+i\beta_r^{k_4}\right)}\right|^2
\right\}^3, \\
& L\to\infty, \\
& 1\leq k_1,k_2,k_3,k_3\leq k_0.
\end{split}
\ee
\end{mydef1}

Next, we give following corollaries from the Theorem, (comp. (b) in the beginning of the section 1.3).

\begin{mydef41}
\be \label{3.2}
\begin{split}
& \prod_{r=1}^{k_2}\left|\frac{\zeta\left(\frac 12+i\alpha_r^{2,k_2}\right)}{\zeta\left(\frac 12+i\beta_r^{k_2}\right)}\right|^2\sim \\
& \sim \frac{\cos^3(\alpha_0^{1,k_1}}{\cos(\alpha_0^{2,k_2})}
\prod_{r=1}^{k_1}\left|\frac{\zeta\left(\frac 12+i\alpha_r^{1,k_1}\right)}{\zeta\left(\frac 12+i\beta_r^{k_1}\right)}\right|^2+
\frac{U^2}{3}\frac{1}{\cos(\alpha_0^{2,k_2})\cos^3U}\times \\
& \times
\left\{
\cos^2(\alpha_0^{3,k_3})\prod_{r=1}^{k_3}\left|\frac{\zeta\left(\frac 12+i\alpha_r^{3,k_3}\right)}{\zeta\left(\frac 12+i\beta_r^{k_3}\right)}\right|^2-
\sin^2(\alpha_0^{4,k_4})\prod_{r=1}^{k_4}
\left|\frac{\zeta\left(\frac 12+i\alpha_r^{4,k_4}\right)}{\zeta\left(\frac 12+i\beta_r^{k_4}\right)}\right|^2
\right\}^3, \\
&  L\to\infty.
\end{split}
\ee
\end{mydef41}

\begin{mydef42}
\be \label{3.3}
\begin{split}
& \prod_{r=1}^{k_3}\left|\frac{\zeta\left(\frac 12+i\alpha_r^{3,k_3}\right)}{\zeta\left(\frac 12+i\beta_r^{k_3}\right)}\right|^2\sim \\
& \sim \frac{\sin^2(\alpha_0^{4,k_4})}{\cos^2(\alpha_0^{3,k_3})}\prod_{r=1}^{k_4}
\left|\frac{\zeta\left(\frac 12+i\alpha_r^{4,k_4}\right)}{\zeta\left(\frac 12+i\beta_r^{k_4}\right)}\right|^2+ \\
& +\left\{
\frac{3}{U^2}\frac{\cos^3U\cos(\alpha_0^{2,k_2})}{\cos^6(\alpha_0^{3,k_3})}
\prod_{r=1}^{k_2}\left|\frac{\zeta\left(\frac 12+i\alpha_r^{2,k_2}\right)}{\zeta\left(\frac 12+i\beta_r^{k_2}\right)}\right|^2-
\right. \\
& \left.
-\frac{3}{U^2}\frac{\cos^3U\cos^3(\alpha_0^{1,k_1})}{\cos^6(\alpha_0^{3,k_3})}
\prod_{r=1}^{k_1}\left|\frac{\zeta\left(\frac 12+i\alpha_r^{1,k_1}\right)}{\zeta\left(\frac 12+i\beta_r^{k_1}\right)}\right|^2
\right\}^{1/3} , \\
& L\to\infty .
\end{split}
\ee
\end{mydef42}

\begin{mydef43}
\be \label{3.4}
\begin{split}
& \prod_{r=1}^{k_4}
\left|\frac{\zeta\left(\frac 12+i\alpha_r^{4,k_4}\right)}{\zeta\left(\frac 12+i\beta_r^{k_4}\right)}\right|^2 \sim \\
& \sim \frac{\cos^2(\alpha_0^{3,k_3})}{\sin^2(\alpha_0^{4,k_4})}\prod_{r=1}^{k_3}
\left|\frac{\zeta\left(\frac 12+i\alpha_r^{3,k_3}\right)}{\zeta\left(\frac 12+i\beta_r^{k_3}\right)}\right|^2- \\
&- \left\{
\frac{3}{U^2}\frac{\cos^3U\cos(\alpha_0^{2,k_2})}{\sin^6(\alpha_0^{4,k_4})}
\prod_{r=1}^{k_2}\left|\frac{\zeta\left(\frac 12+i\alpha_r^{2,k_2}\right)}{\zeta\left(\frac 12+i\beta_r^{k_2}\right)}\right|^2-
\right. \\
& \left.
-\frac{3}{U^2}\frac{\cos^3U\cos^3(\alpha_0^{1,k_1})}{\sin^6(\alpha_0^{4,k_4})}
\prod_{r=1}^{k_1}\left|\frac{\zeta\left(\frac 12+i\alpha_r^{1,k_1}\right)}{\zeta\left(\frac 12+i\beta_r^{k_1}\right)}\right|^2
\right\}^{1/3} , \\
& L\to\infty .
\end{split}
\ee
\end{mydef43}

\subsection{}

Let us remind the following correspondences (see (\ref{2.1}), (\ref{2.2}); (\ref{2.3}), (\ref{2.4}); (\ref{2.9})):
\be \label{3.5}
\begin{split}
& f_1(t)=\cos^3t \rightarrow \prod_{r=1}^{k_1}\left|\frac{\zeta\left(\frac 12+i\alpha_r^{1,k_1}\right)}{\zeta\left(\frac 12+i\beta_r^{k_1}\right)}\right|^2=(k_1,f_1), \\
& f_2(t)=\cos t \rightarrow \prod_{r=1}^{k_2}\left|\frac{\zeta\left(\frac 12+i\alpha_r^{2,k_2}\right)}{\zeta\left(\frac 12+i\beta_r^{k_2}\right)}\right|^2 = (k_2,f_2), \\
& f_3(t)=\cos^2t \rightarrow \prod_{r=1}^{k_3}\left|\frac{\zeta\left(\frac 12+i\alpha_r^{3,k_3}\right)}{\zeta\left(\frac 12+i\beta_r^{k_3}\right)}\right|^2 = (k_3,f_3), \\
& f_4(t)=\sin^2t \rightarrow \prod_{r=1}^{k_4}\left|\frac{\zeta\left(\frac 12+i\alpha_r^{4,k_4}\right)}{\zeta\left(\frac 12+i\beta_r^{k_4}\right)}\right|^2 = (k_4,f_4), \\
& 1\leq k_1,k_2,k_3,k_4\leq k_0.
\end{split}
\ee

\begin{remark}
First, we see that the formulae (\ref{3.1}) -- (\ref{3.4}) define the set of nonlinear interactions of the $(k_0)^4$ oscillating systems in (\ref{3.5}). Namely, this set
contains $4(k_0)^4$ elements of different type. Now, if we use our short notions in (\ref{3.5}) we may write down the following diagram
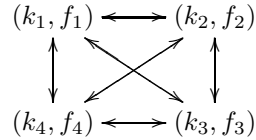
\begin{displaymath}
    \xymatrix{
        (k_1,f_1) \ar[d] \ar[dr] \ar[r]  & (k_2,f_2) \ar[d] \ar[l] \ar[dl] \\
        (k_4,f_4) \ar[u] \ar[ur] \ar[r]   & (k_3,f_3) \ar[u] \ar[l] \ar[ul] }
\end{displaymath}
\end{remark}

\section{Short functional algebra}

We shall call each of the following type of functional combinations
\be \label{4.1}
\cos^2(\alpha_0^{3,k_3})\prod_1^{k_3}-\sin^2(\alpha_0^{4,k_4})\prod_1^{k_4} ,
\ee
\be \label{4.2}
\frac{3}{U^2}\cos^3U\cos(\alpha_0^{2,k_2})\prod_1^{k_2}-\frac{3}{U^2}\cos^3U\cos^3(\alpha_0^{1,k_1})\prod_1^{k_1}
\ee
of two $\zeta$-oscillating systems as the two-parts form. Since
\bdis
1\leq k_1,k_2,k_3,k_4\leq k_0
\edis
then the set of formulae (\ref{4.1}) (just as similar set (\ref{4.2})) contains
\bdis
(k_0)^2
\edis
two-part forms of different type. In this direction we have (see (\ref{3.1})) the following

\begin{mydef44}
\be \label{4.3}
\begin{split}
& \left\{
\cos^2(\alpha_0^{3,k_3})\prod_{r=1}^{k_3}\left|\frac{\zeta\left(\frac 12+i\alpha_r^{3,k_3}\right)}{\zeta\left(\frac 12+i\beta_r^{k_3}\right)}\right|^2 -
\right. \\
& \left. -
\sin^2(\alpha_0^{4,k_4})\prod_{r=1}^{k_4}\left|\frac{\zeta\left(\frac 12+i\alpha_r^{4,k_4}\right)}{\zeta\left(\frac 12+i\beta_r^{k_4}\right)}\right|^2
\right\}^3 \sim \\
& \sim
\frac{3}{U^2}\cos^3U\cos(\alpha_0^{2,k_2})\prod_{r=1}^{k_2}\left|\frac{\zeta\left(\frac 12+i\alpha_r^{2,k_2}\right)}{\zeta\left(\frac 12+i\beta_r^{k_2}\right)}\right|^2- \\
& -\frac{3}{U^2}\cos^3U\cos^3(\alpha_0^{1,k_1})\prod_{r=1}^{k_1}\left|\frac{\zeta\left(\frac 12+i\alpha_r^{1,k_1}\right)}{\zeta\left(\frac 12+i\beta_r^{k_1}\right)}\right|^2, \\
& L\to\infty .
\end{split}
\ee
\end{mydef44}

\begin{remark}
The following property is expressed by the formula (\ref{4.3}): the cube of two-parts form (\ref{4.1}) is asymptotically equal to the two-part form (\ref{4.2}). Of course, we have also the
following formula
\bdis
\sqrt[3]{(4.2)}\sim (4.1).
\edis
\end{remark}

Next, it is true (see (\ref{4.3})) that for every fixed pair $(k_3,k_4)$ we have set of $(k_0)^2$ asymptotic formulae of different types for the cube of corresponding two-parts form (\ref{4.1}).
For example, in the case $k_0=10^3$ we have $10^6$ of these formulae.

\begin{remark}
Consequently, we may regard the formula (\ref{4.3}) as a kind of simplification of the well-known school-formula
\be \label{4.4}
(a-b)^3=a^3-3a^2b+2ab^2-b^3
\ee
in this short functional algebra generated by the formula (\ref{3.1}). Namely, the right-hand side of (\ref{4.3}) contains two-parts form (i.e. the two members only in comparison with
(\ref{4.4})).
\end{remark}

\section{The iteration formula as a constraint on the corresponding vector-valued functions generated by the operator $\hat{H}$}

Let us remind that we have defined (see \cite{8}, Definition 2 and Definition 5; \cite{9}, Definition) new type of vector-valued operator $\hat{H}$ as follows:
\bdis
\begin{split}
& \forall\- f(t)\in \tilde{C}_0[T,T+U] \rightarrow \hat{H}f(t)=\\
& = (\alpha_0,\alpha_1,\dots,\alpha_k,\beta_1,\dots,\beta_k),\ k=1,\dots,k_0,
\end{split}
\edis
for every fixed $k$.

Consequently, we have the following set of $k_0$ vector-valued functions
\be \label{5.1}
(\alpha_0,\alpha_1,\dots,\alpha_k,\beta_1,\dots,\beta_k),\ k=1,\dots,k_0 .
\ee

\begin{remark}
That is, we may say that we have defined the following matrix-valued operator $\hat{H}$:
\bdis
f(t)\overset{\hat{H}}{\rightarrow}
\begin{pmatrix}
\alpha_0 & \alpha_1 & \beta_ 1 & 0 & \dots \\
\alpha_0 & \alpha_1 & \alpha_2 & \beta_1 & \beta _2 & 0 & \dots \\
\alpha_0 & \alpha_1 & \alpha_2 & \alpha_3 & \beta_1 & \beta _2 & \beta_3 & 0 & \dots \\
\vdots \\
\alpha_0 & \alpha_1 & \alpha_2 & \alpha_3 & \dots & \alpha_{k_0} & \beta_1 & \beta_2 & \dots & \beta_{k_0}
\end{pmatrix}_{k_0 \times (2k_0+1)} .
\edis
\end{remark}

Next, it is true that every interaction formula:
\begin{itemize}
\item[(a)] contains some set of $\zeta$-oscillating systems,
\item[(b)] every $\zeta$-oscillating system from that set contains the components of corresponding vector-valued function of type (\ref{5.1}).
\end{itemize}

\begin{remark}
Consequently, we may understand every interaction formula (\ref{3.1}), for example, as the constraint on the set of corresponding vector-valued
functions of type (\ref{5.1}) which are contained in this formula.
\end{remark}


\end{document}